\documentclass[switch]{article} 

\usepackage[utf8]{inputenc} 
\usepackage[T1]{fontenc}    
\usepackage{hyperref}

\usepackage{enumitem}
\usepackage{url}            
\usepackage{booktabs}       
\usepackage{bm}
\usepackage{amsfonts}       
\usepackage{amsmath}
\usepackage{amssymb}
\usepackage{nicefrac}       
\usepackage{microtype}
\usepackage{cleveref}
\usepackage{xcolor}
\usepackage{algpseudocode}
\usepackage{algorithm}
\usepackage{color,soul}
\usepackage{caption}
\usepackage{lipsum}
\usepackage{dsfont}
\usepackage{breqn}
\usepackage{comment}
\usepackage{multicol}
\usepackage{stfloats}

\usepackage{amsmath,amssymb}
\usepackage{fancyhdr}       
\usepackage{graphicx}       
\graphicspath{{media/}}     

\makeatletter
\newcommand{\printfnsymbol}[1]{%
  \textsuperscript{\@fnsymbol{#1}}%
}

\newcommand*\bigcdot{\mathpalette\bigcdot@{.5}}
\newcommand*\bigcdot@[2]{\mathbin{\vcenter{\hbox{\scalebox{#2}{$\m@th#1\bullet$}}}}}

\NewDocumentCommand{\ARef}{ s s m }{%
    \IfBooleanTF{#2}{}{%
        \cref{#3}%
    }%
    \IfBooleanT{#1}{%
        \IfBooleanF{#2}{%
            , %
        }%
        line~\ref{#3}%
    }%
}

\makeatother

\usepackage{arxiv}
\usepackage{graphicx} 

\title{PhoTOS: Topology Optimization of Photonic Components using a Shape Library}

\author{
\and
Rahul Kumar Padhy \\
Department of Mechanical Engineering \\
University of Wisconsin-Madison \\
Madison, WI, USA \\
\texttt{rkpadhy@wisc.edu} \\
\and
Aaditya Chandrasekhar \\
Department of Mechanical Engineering \\
Northwestern University \\
Evanston, IL, USA \\
\texttt{cs.aaditya@gmail.com} \\
}

\begin{document}


\maketitle

\begin{abstract}
Topology Optimization (TO) holds the promise of designing next-generation compact and efficient photonic components. However, ensuring the optimized designs comply with fabrication constraints imposed by semiconductor foundries remains a challenge. This work presents a TO framework that guarantees designs satisfy fabrication criteria, particularly minimum feature size and separation. Leveraging recent advancements in machine learning and feature mapping methods, our approach constructs components by transforming shapes from a predefined library, simplifying constraint enforcement. Specifically, we introduce a Convo-implicit Variational Autoencoder to encode the discrete shape library into a differentiable space, enabling gradient-based optimization. The efficacy of our framework is demonstrated through the design of several common photonic components.
\end{abstract}


 \vspace{1cm}
 

\section{Introduction}
\label{sec:intro}

\begin{figure}[h]
 	\begin{center}

    \includegraphics[scale=0.26,trim={0 0 0 0},clip]{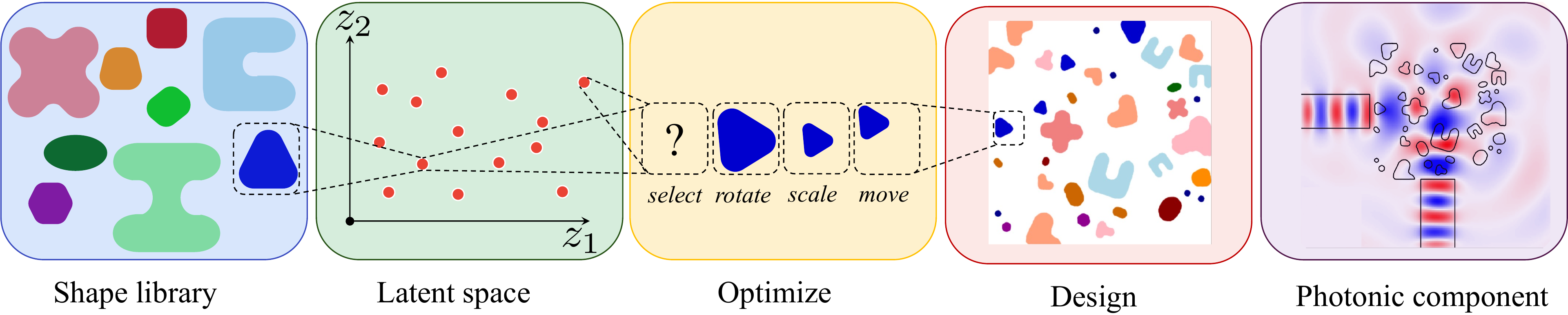}
 		\caption{Graphical abstract: Given a predefined library of shapes, a Convo-implicit Variational Autoencoder (VAE) is trained to encode them in a differentiable latent space. Shape instances are then selected from this latent space and subjected to rotation, scaling, and translation to populate the design space, yielding an optimized photonic component.}
        \label{fig:graphical_abstract}
	\end{center}
 \end{figure}

Topology optimization (TO) are a class of methods used to optimize material distribution within a design domain, achieving optimal performance under given constraints. These methods efficiently navigate complex design spaces that challenge conventional, intuition-based approaches. While traditionally applied in structural mechanics, TO has seen increasing adoption in photonics \cite{christiansen2021inverseTutorialPhotonics}, where advances in integrated photonic device design are critical for progress in high-speed communication \cite{marpaung2019integrated}, quantum computing \cite{arrazola2021quantum}, and machine learning accelerators \cite{wetzstein2020inference}.

\begin{figure}
 	\begin{center}
		\includegraphics[scale=0.5,trim={0 0 0 0},clip]{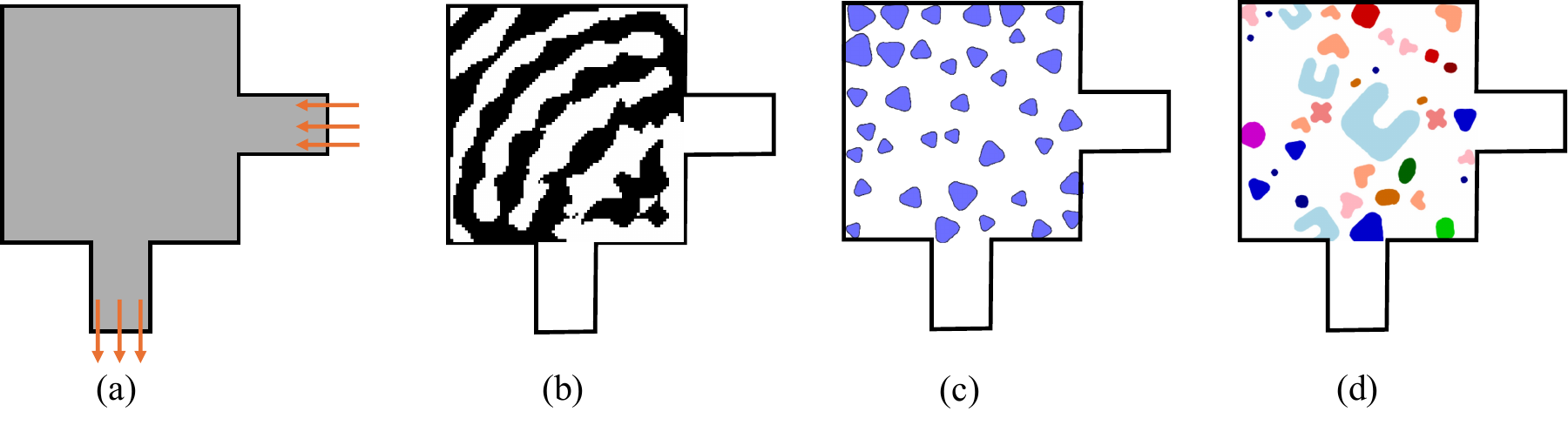}
 		\caption{(a) Photonic device design domain and boundary conditions.
                        (b) Density-based topology optimization. 
                        (c) Topology optimization using single shape-based feature mapping.
                        (d) Topology optimization using multiple generic shapes.}
        \label{fig:TO_methods}
	\end{center}
 \end{figure}

While various TO methods have been proposed, density-based approaches have been widely adopted \cite{sigmund2013topology}. Here, the design is parameterized by a set of pixels, and the optimal material is assigned to each pixel \cite{jenkins2022general, chen2020design}. This results in organic, free-form designs. Consider, for example, the design of a waveguide bend (\cref{fig:TO_methods}(a)), where light from the input port (right) is to be directed to the output port (bottom) with minimal back reflection and maximum transmission. \cref{fig:TO_methods}(b) illustrates a typical density-based design. While offering significant design freedom, the resulting designs can be difficult to interpret \cite{norato2018SuperShapes}, manufacture, and constrain to fabrication criteria \cite{schubert2022inverse, chen2020design, piggott2015inverse}.

Conversely, alternative techniques have been proposed that represent the design using simpler geometric shapes such as bars \cite{norato2015gpto}, plates \cite{zhang2016PlateGeomProj}, and polygons \cite{chandrasekhar2023polyto}. These methods, collectively referred to as feature mapping methods \cite{wein2020reviewFeatureMapping} parameterize shapes by high-level descriptors such as width, radius, and angle. Varying these parameters generates various dimensions of the shape that are projected onto the design domain. For instance \cref{fig:TO_methods}(c) showcases the waveguide bend composed of triangles. As evident, these techniques yield more interpretable designs that are easier to constrain. However, while simplifying the enforcement of manufacturing constraints, they significantly limit design freedom \cite{norato2018SuperShapes}.

This work proposes an optimization framework that strikes a balance between density-based and feature-mapping methods. Our approach offers increased design freedom \cite{choi2021curvilinear} while retaining the interpretability and ease of constraint enforcement characteristic of feature-mapping techniques. Specifically, we make the following contributions:

\begin{enumerate}
    \item We extend feature-mapping methods, traditionally limited to a single shape, to accommodate multiple shapes (\cref{sec:shape_shapeLibrary}).
    \item We introduce a Convo-implicit Variational Autoencoder (VAE) \cite{park2019deepsdf} to transform a discrete library of shapes into a continuous and differentiable form, enabling gradient based optimization (\cref{sec:method_libraryRep_convoImpNN}).
    \item We present an optimization framework that produce fabricable, optimized photonic components by selecting and transforming shapes from the library (\cref{sec:method_optimization}).
    \item We apply our framework to design waveguide bends and mode converters, demonstrating its ability to reliably produce photonic components that meet both performance and fabrication criteria (\cref{sec:results}).
\end{enumerate}

\section{Proposed Method}
\label{sec:method}

\subsection{Overview}
\label{sec:method_overview}

\begin{figure}
 	\begin{center}
		\includegraphics[scale=0.45,trim={0 0 0 0},clip]{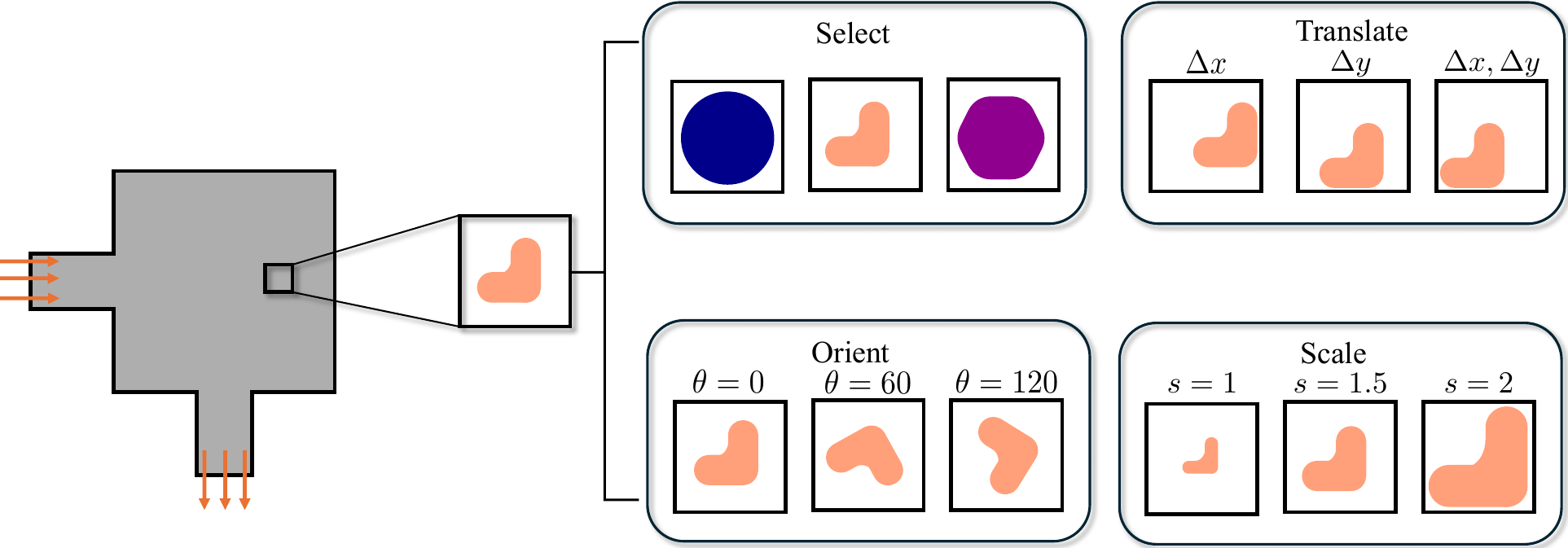}
 		\caption{Shape instances selected from the library are translated, oriented, and scaled onto the design domain.}
            \label{fig:variation_transform_parameters}
	\end{center}
 \end{figure}

In this study, we focus on the TO of photonic components. We begin by assuming that the performance criteria; typically the allowable insertion, reflection and cross-talk loss have been prescribed. Further, we assume that the fabrication criteria: minimum feature size (MFS) and minimum separation distance (MSD) (\cref{fig:min_feature_size}) have also been specified. Finally, we assume that a library of ($n_L$) shapes has been prescribed (\cref{fig:feature_library}). The objective then is to find through gradient based optimization, an optimal configuration (selection, translation, orientation, and scaling) (see \cref{fig:variation_transform_parameters}) of shape instances from the library onto the design space such that the optimized design meets the required performance and fabrication criteria.

\subsection{Shape Library}
\label{sec:shape_shapeLibrary}

\begin{figure}
 	\begin{center}
		\includegraphics[scale=0.5,trim={0 0 0 0},clip]{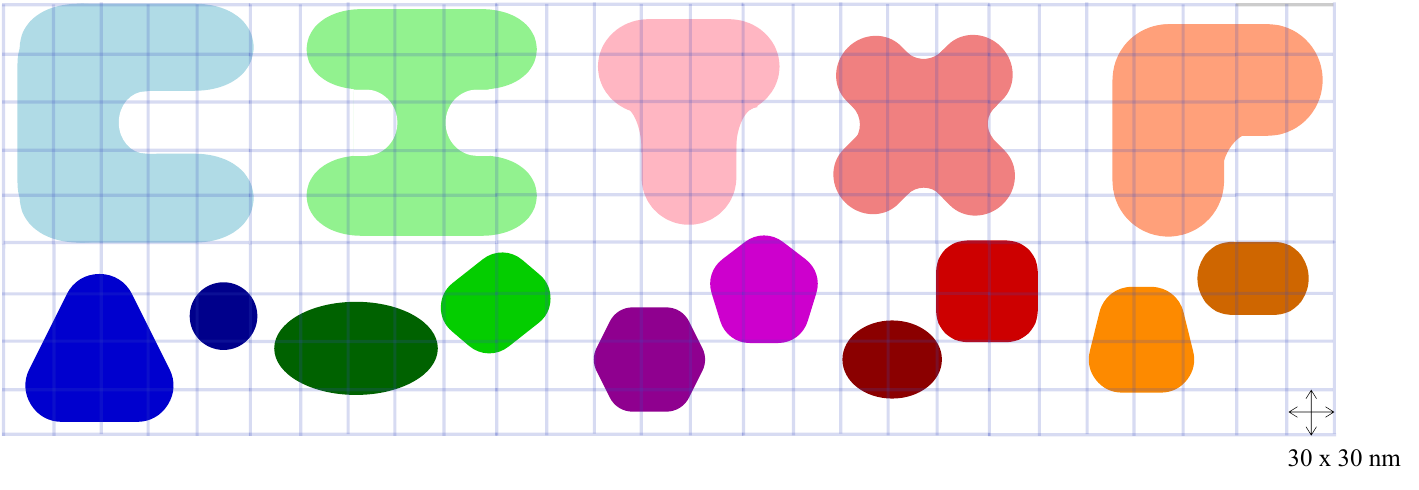}
 		\caption{A library of $n_L \; (=15)$ prescribed shapes.}
        \label{fig:feature_library}
	\end{center}
 \end{figure}
 
\begin{figure}
 	\begin{center}
        \includegraphics[scale=1.2,trim={0 0 0 0},clip]{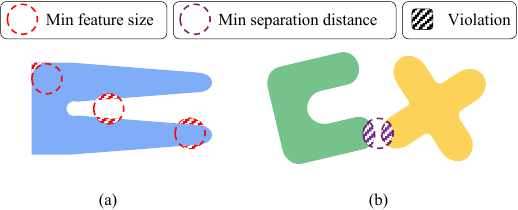}
 		\caption{To ensure fabricability, we impose: (a) minimum feature size constraint and (b) minimum separation distance constraints. }
        \label{fig:min_feature_size}
	\end{center}
 \end{figure}

We compose our library (\cref{fig:feature_library}) of generic shapes frequently employed in photomask fabrication \cite{puthankovilakam2017unified} and other common shapes. The shapes are chosen to meet an MFS (\cref{fig:min_feature_size}(a)) of 40 nm \cite{chen2024PhotonicsValidation}. Importantly, we observe that these shapes are not derived from a single parametric family; rather, the library includes a variety of $\textit{discrete}$ shapes without common defining parameters.

Recall that gradient-based optimization relies on the ability to $\textit{continuously}$ vary its optimization parameters. In other words, the optimizer must be able to continuously vary between the shapes in the library. Thus, to facilitate the optimal choice of shapes by the optimizer, we transform our $\textit{discrete}$ shape library into a continuous and differentiable representation. In particular, a continuous representation of the shape library entails two key components:
\begin{enumerate}[label=(\alph*)]
    \item a continuous representation of each individual shape within the library (\cref{sec:method_libraryRep_sdf}).
    \item a continuous representation of the collective library itself (\cref{sec:method_libraryRep_convoImpNN}).
\end{enumerate}

\subsubsection{Shape Representation}
\label{sec:method_libraryRep_sdf}

\begin{figure}
 	\begin{center}
		\includegraphics[scale=0.75,trim={0 0 0 0},clip]{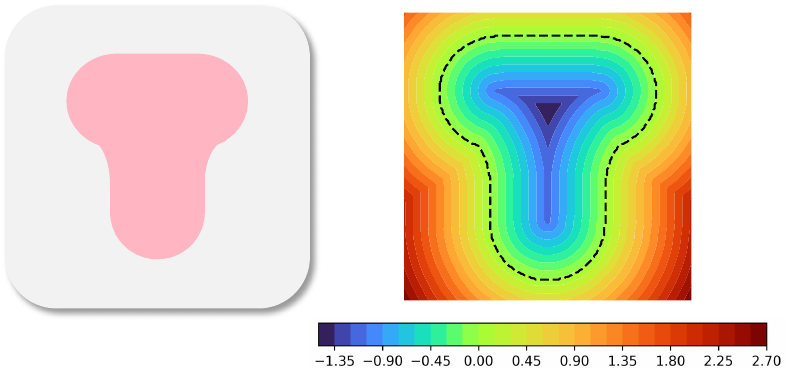}
 		\caption{A shape from the library and its signed distance field (SDF). The $0^{th}$ contour of the SDF is shown in dotted lines.}
        \label{fig:shape_and_sdf}
	\end{center}
 \end{figure}
 
To achieve continuous shape representation, we represent our shapes using signed distance fields (SDFs). SDFs ($\phi$) are $C^0$ continuous functions, with the value at any point defined as the shortest distance to the shape's boundary (\cref{fig:shape_and_sdf}). We compute the SDFs for the shapes in our library over the bounding box ($[\underline{\bm{x}}_{b}, \overline{\bm{x}}_{b}]$).  In addition to their continuous nature, the implicit definition of distances from the boundary, facilitates the computation of the MSD constraint (see \cref{fig:min_feature_size}(b), \cref{sec:method_optimization}).

\subsubsection{Library Representation}
\label{sec:method_libraryRep_convoImpNN}

To achieve continuous library representation, we employ a VAE. VAEs 
are a particular construct of neural networks, that among other attributes, convert discrete data into a continuous and differentiable representation. This facilitates gradient-based optimization \cite{kingma2019introduction, kingma2013auto, doersch2016tutorial}.

Specifically, the shape SDFs from the previous section serve as both the input and output data for the VAE. We propose a VAE architecture with the following key attributes:

\begin{enumerate}[label=(\alph*)] 

\item A continuous and differentiable latent space that represents the shape library, facilitating gradient-based optimization.

\item An implicit representation of the output shape SDFs as a function of spatial coordinates. This enables the querying of SDF values at arbitrary locations and resolutions, a crucial requirement for our optimization formulation.

\end{enumerate}

\begin{figure}
 	\begin{center}
            \includegraphics[scale=0.65,trim={0 0 0 0},clip]{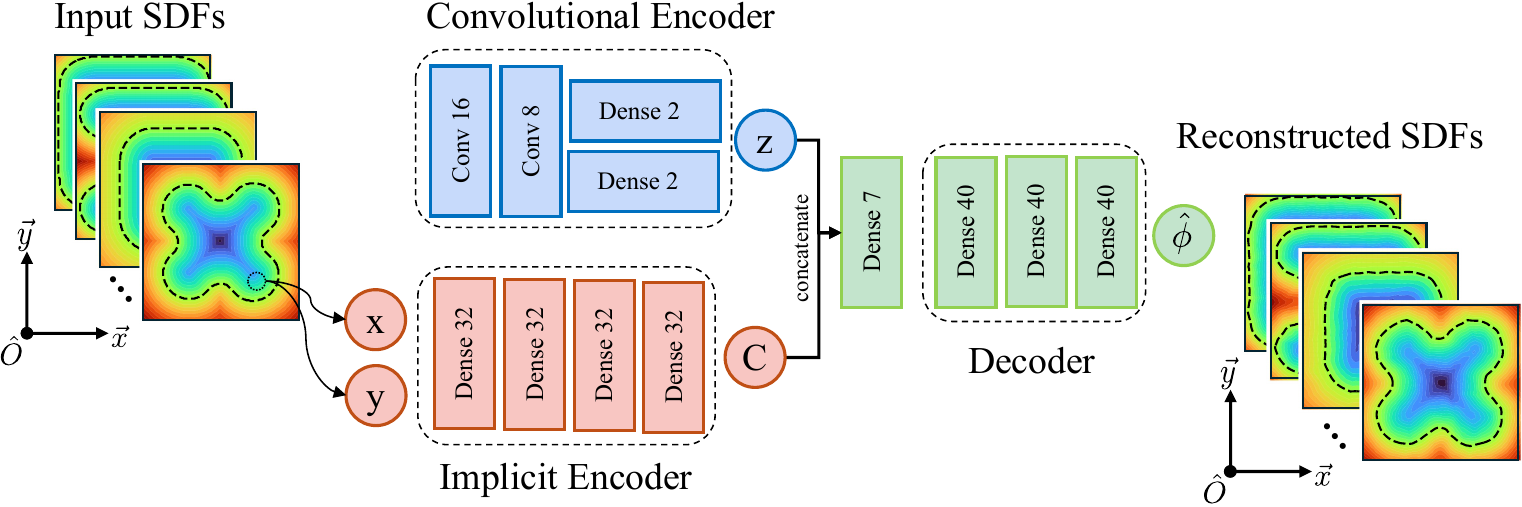}
 		\caption{Architecture of the convo-implicit VAE.}
        \label{fig:photos_NN}
	\end{center}
 \end{figure}

To achieve these attributes, we propose a convo-implicit VAE architecture. This architecture employs a convolutional encoder in conjunction with an implicit encoder \cite{park2019deepsdf, sitzmann2020SIREN, tancik2020fourier} and decoder. The proposed Convo-Implicit VAE architecture (\cref{fig:photos_NN}) comprises the following components:

\begin{enumerate}
    \item Firstly, the shape SDFs $( \mathcal{I} = \{ \phi^{(1)}, \phi^{(2)}, \ldots, \phi^{(L)} \})$ of size $\{n_L \times n_x \times n_x \}$ are propagated through the convolutional encoder. Here, $n_x$ denotes the resolution of the SDF images in the training set. This results in an output latent space $z$ of size $\{n_L \times n_z\}$. In particular, we have $n_z = 2$.

    \item The coordinates of the SDFs ($\bm{x}$) of size $\{n_p \times 2\}$ are propagated through an implicit encoder \cite{sitzmann2020SIREN}. Here, $n_p = n_x^2$ is the total number of pixels in the training images. This results in an output of projected coordinates $C$ of size $\{n_p \times n_c\}$.
    
    \item The latent space $z$ and the projected coordinates $C$ are concatenated. This combined tensor of size $\{n_L \times n_p \times (n_z + n_c)\}$ is propagated through the decoder, producing the reconstructed feature SDFs $(\hat{\mathcal{I}} = \{\hat{\phi}^{(1)}, \hat{\phi}^{(2)}, \ldots, \hat{\phi}^{(L)}\})$ of size $\{n_L \times n_p\}$.

\end{enumerate}

\begin{figure}
 	\begin{center}
		\includegraphics[scale=0.3,trim={0 0 0 0},clip]{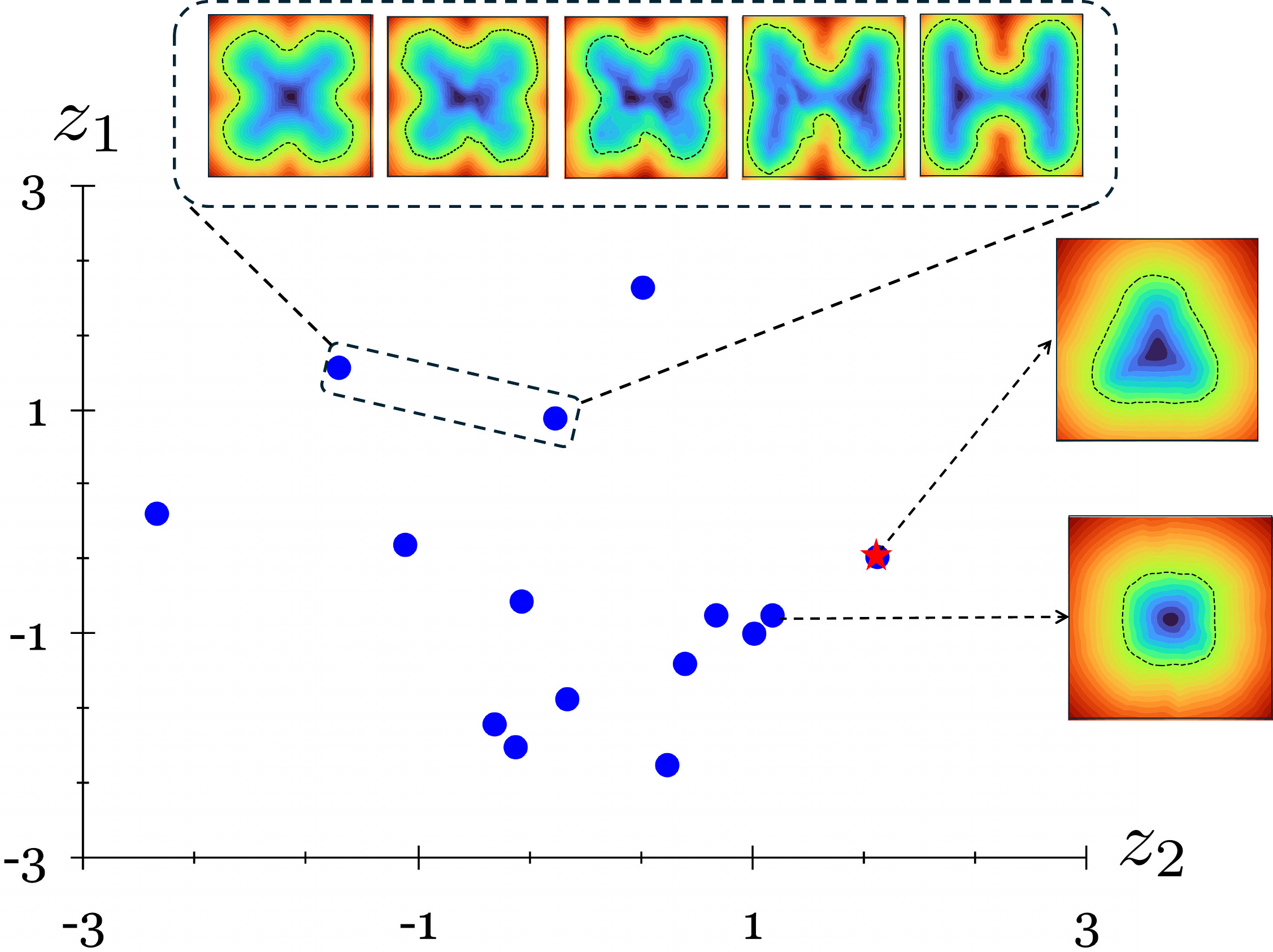}
 		\caption{Continuous latent space obtained from the trained Conv-implicit VAE. (\textcolor{red}{$\bigstar$}) represents the coordinates $[1.7, -0.3]$, corresponding to a rounded triangle.}
        \label{fig:latent_space}
	\end{center}
 \end{figure}
 
We train the Convo-Implicit VAE to minimize the discrepancy between the reconstructed and input SDFs by reducing their mean squared error \cite{rautela2022towards}. Additionally, the latent space (\cref{fig:latent_space}) is constrained to resemble a standard normal distribution $\bm{z} \sim \mathcal{N}(0, 1)$ through a KL divergence loss term \cite{kingma2019introduction}. The Convo-Implicit VAE's net loss can be expressed as:

\begin{equation}
    L_v = ||\mathcal{I} - \hat{\mathcal{I}} ||_2^2 \;  + \; \beta \text{KL}(\bm{z} || \mathcal{N})
    \label{eq:vae_training_loss}
\end{equation}
Where the parameter $\beta$ $(=5 \times 10^{-8})$ is the relative weight of the KL divergence loss term. Convergence was achieved after training for $4\times 10^{4}$ epochs with a learning rate of $5\times 10^{-4}$ using the Adam optimizer.

After training, we discard the convolutional encoder. The retained implicit encoder and decoder can be queried with spatial $\bm{x}$ and latent $\bm{z}$ coordinates respectively to obtain the SDF $\hat{\phi}(\bm{x}, \bm{z})$. For instance querying at $\bm{z} = [1.7, -0.3]$ and $\bm{x} \in [\underline{\bm{x}}_{b}, \overline{\bm{x}}_{b} ]$, we obtain the SDF of the rounded triangle (\cref{fig:latent_space}).

\subsubsection{Approximating the SDF}
\label{sec:method_libraryRep_approxSDF}

The retained implicit encoder and decoder allow us to query shape SDFs at any spatial coordinate, including those outside the training bounding box (\cref{sec:shape_shapeLibrary}). However, this presents two challenges:
\begin{enumerate}
    \item Querying the VAE for each shape instance over a large number of spatial coordinates can be computationally expensive.
    \item The VAE must extrapolate SDF values outside the training bounding box, leading to inaccurate results.
\end{enumerate}

\begin{figure}
 	\begin{center}
		\includegraphics[scale=0.5,trim={0 0 0 0},clip]{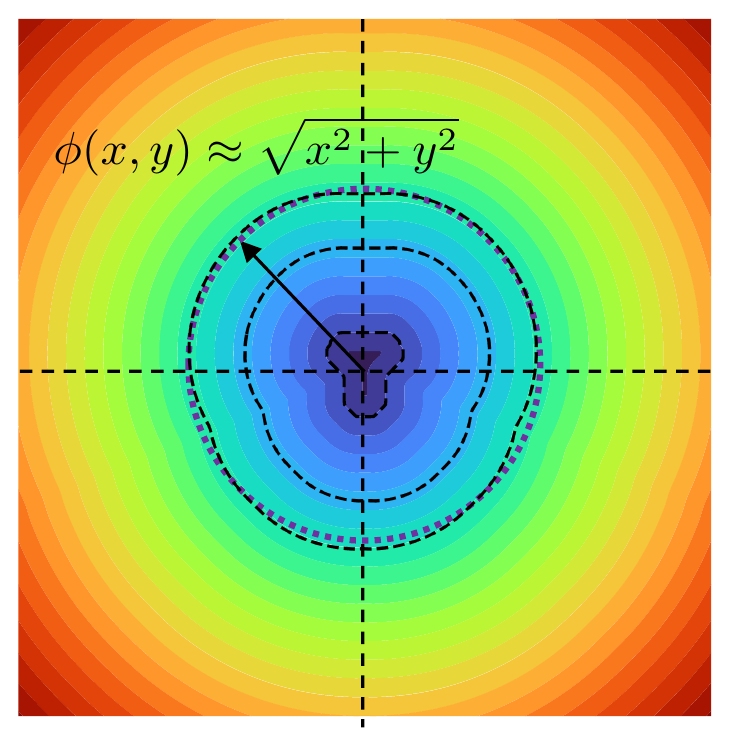}
 		\caption{SDF at large distance distance from an object can be approximated to a circle.}
        \label{fig:far_sdf}
	\end{center}
 \end{figure}

To address these challenges, we propose approximating the shape SDF outside its bounding box. Observe that at large distances, the SDF of any shape approaches the Euclidean distance from the origin (Figure \ref{fig:far_sdf}). Therefore, we can approximate the SDF as:

\begin{equation}
     \hat{\phi}(\bm{x}) = \begin{cases}
                        \hat{\phi}(\bm{x}) \quad \bm{x} \in [\underline{\bm{x}}_{b}, \overline{\bm{x}}_{b} ] \\
                        \approx ||\bm{x}||_2 \quad \text{Otherwise}
                    \end{cases}
    \label{eq:approx_sdf}
\end{equation}

Importantly, note that since we retain the SDF values within the training bounding box, this approximation does not compromise the accuracy of our optimization.

\subsubsection{Geometric Projection}
\label{sec:method_projection}

\begin{figure}
 	\begin{center}
		\includegraphics[scale=0.75,trim={0 0 0 0},clip]{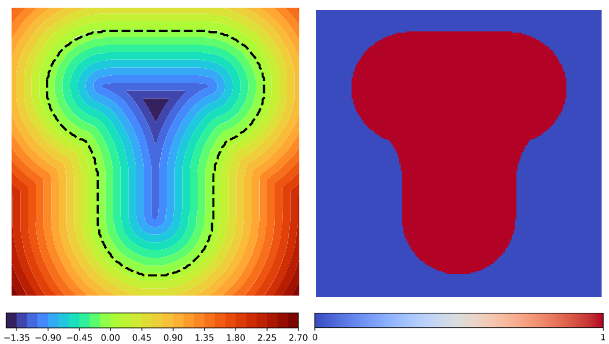}
 		\caption{The SDF of a shape and its projected density field.}
        \label{fig:sdf}
	\end{center}
 \end{figure}
 
The central idea of geometry projection is to map shapes, as described by their SDFs, onto density fields (\cref{fig:sdf}). This mapping subsequently enables us to simulate the design. The density field, $\rho(x)$, is defined such that $\rho(\bm{x}) = 0$ denotes the region outside the shape and $\rho(\bm{x}) = 1$ denotes the region inside the shape. Given $\hat{\phi}^{(i)}(\bm{x})$ as the SDF at a point $\bm{x}$ of shape $i$, the corresponding density $\rho^{(i)}(\bm{x})$, is obtained as:

\begin{equation}
    \rho^{(i)}(\bm{x}) = \sigma_{\beta}(\hat{\phi}(\bm{x}^{(i)})) = \frac{1}{1 + e^{-\beta \hat{\phi}(\bm{x}^{(i)})}}
    \label{eq:density_projection}
\end{equation}

 Where $\sigma_{\beta}$ is the sigmoid function with a projection sharpness of $\beta$. Observe that this mapping is differentiable, thus enabling gradient-based optimization. Once the density fields of all shape instances are computed, the overall design density $\hat{\rho}(\bm{x})$ is obtained by taking the union using a $p^{th}$-norm formulation:

\begin{equation}
    \hat{\rho}(\bm{x}) = \left( \sum\limits_{i=1}^{n_F} \left( \rho^{(i)}(\bm{x} \big)\right)^p \right)^{\frac{1}{p}}
    \label{eq:density_union}
\end{equation}

Finally, the design domain is populated with Silicon where the density equals 1, and Silicon dioxide where the density equals 0. The permittivity at a point $\bm{x}$ can then be expressed as:

\begin{equation}
    \varepsilon(\bm{x}) = \varepsilon^{(ox)} + (\varepsilon^{(si)} - \varepsilon^{(ox)})\hat{\rho}(\bm{x})
    \label{eq:linear_material}
\end{equation}

Where $ \varepsilon^{(ox)} (= 2.25)$ and $\varepsilon^{(si)} (= 12.25)$ are the permittivities of Silicon oxide and Silicon respectively.

\subsection{Shape Transformation}
\label{sec:method_Transformation}

Recall that our objective is to achieve an optimal $\textit{configuration}$ of shapes from the library within the design space (\cref{fig:variation_transform_parameters}). Further, recall that each shape is defined in its local (library) frame of reference (\cref{fig:photos_NN}). We achieve the configurations through affine transformations from the local to the global (design) frame of reference. In particular, we define a translation ($\bar{x}^{(i)}, \bar{y}^{(i)}$), orientation ($\theta^{(i)}$), and scaling ($s^{(i)}$) operation for each shape instances $i = 1,\ldots, n_F$. With ${x}_e ,{y}_e$ being coordinates in the design frame of reference (typically, center coordinates of a simulation mesh elements), we obtain the transformed coordinates $ \tilde{x}_e^{(i)}, \tilde{y}_e^{(i)}$ as:

\begin{equation}
    \begin{pmatrix}
    \tilde{x}_e^{(i)} \\ \tilde{y}_e^{(i)}
    \end{pmatrix} = \frac{1}{s^{(i)}}
    \bm{R}(\theta^{(i)})
    \begin{pmatrix}
    {x_t}_e^{(i)} \\ {y_t}_e^{(i)}
    \end{pmatrix} 
    \label{eq:affine_transform}
\end{equation}

where $ {x_t}_e^{(i)} , {y_t}_e^{(i)} $ are the coordinate after translation:
\begin{equation}
    \begin{pmatrix}
    {x_t}_e^{(i)} \\ {y_t}_e^{(i)}
    \end{pmatrix} = 
    \begin{pmatrix}
    x_e - \bar{x}^{(i)} \\ y_e - \bar{y}^{(i)}
    \end{pmatrix} ,
    \label{eq:translation_transform}
\end{equation}

and \( \bm{R}(\theta^{(i)}) \) is the orientation matrix:

\begin{equation}
    \bm{R}(\theta^{(i)}) = 
    \begin{bmatrix}
    \cos \theta^{(i)} & -\sin \theta^{(i)} \\
    \sin \theta^{(i)} & \cos \theta^{(i)} \\
    \end{bmatrix}
    \label{eq:rotation_transform}
\end{equation}

\subsection{Optimization}
\label{sec:method_optimization}

Having established a differentiable representation of the shape library and the necessary methods for its transformation within the design domain, we now outline the key components of the design optimization framework.

\textbf{Design Variables :} The optimization process encompasses the selection (as determined by latent coordinate $(z_1^{(i)}, z_2^{(i)})$), translation $(\overline{x}^{(i)}, \overline{y}^{(i)})$, orientation $(\theta^{(i)})$, and scaling $(s^{(i)})$ of each shape instance $i = 1,\ldots, n_F$. Collectively, this forms our design variables $\bm{\overline{\Lambda}} = \{\bm{\Lambda}^{(1)}, \bm{\Lambda}^{(2)}, \ldots, \bm{\Lambda}^{(n_F)}\}$ where  $\bm{\Lambda}^{(i)} = \{z_1^{(i)}, z_2^{(i)}, \bar{x}^{(i)}, \bar{y}^{(i)}, \theta^{(i)}, s^{(i)}\} $.

\textbf{Simulation :} Upon computing the permittivities (\cref{eq:linear_material}) at the center of mesh elements, we simulate the optical response of the component by employing a Frequency Domain Finite Difference (FDFD) solver.  The electric field, $\bm{E}$, is determined by solving the linear system:

\begin{equation}
    \bm{K}_{\omega} (\bm{\varepsilon}) \bm{E} = \bm{J}
    \label{eq:FDFD_simulation}
\end{equation},

where $\bm{K}_{\omega}$ represents the system matrix at frequency $\omega$ and $\bm{J}$ denotes the source(s). Specifically, we utilize Ceviche \cite{hughes2019forward}, an open-source FDFD electro-magnetic simulator. The resulting electric field is then used to compute the scattering spectra, which serves as the metric for evaluating the objective function.

\textbf{Objective:} We aim to achieve photonic devices that meet performance specifications, specifically, desired levels of insertion and reflection losses. Following the approach in \cite{schubert2022inverse}, we express these specifications by defining cutoff values $\bm{S}^*$ for the magnitudes of the scattering spectra at the monitored wavelengths at the input and output ports. With $\bm{S}(\bar{\bm{\Lambda}})$ representing the scattering spectra obtained with design parameters $\bar{\bm{\Lambda}}$, we can express our objective as:

\begin{equation}
    L(\bar{\bm{\Lambda}}) = \left\| \Psi^+\left(\frac{|\bm{S}(\bar{\bm{\Lambda}})|^2 - |\bm{S}^*|^2}{\bm{w}} \right) \right \|_{2}^2
    \label{eq:objective}
\end{equation}

where $\Psi^+(\cdot)$ denotes the softplus function, and $\bm{w}$ represents the relative weights for each entry in the scattering spectra. For a detailed discussion, we refer the reader to \cite{schubert2022inverse}.

\textbf{Minimum Separation Distance Constraint :}
We impose a minimum spacing constraint to promote fabricability. Following \cite{kang2016structural} \cite{deng2020cvxnet}, we offset the feature instances by an amount equal to half of the MSD (i.e., $\delta^* = MSD/2$). Then, with $\chi^+(\cdot)$ representing the ReLU function and $n_e$ representing the number of elements on the mesh, we can express the minimum separation distance constraint as:

    \begin{equation}
        g_s(\bm{\overline{\Lambda}}) \equiv \frac{1}{n_e} \left( \sum\limits_{e=1}^{n_e} \chi^+\left( \sum\limits_{i=1}^{n_F} \sigma_{\beta}({\phi^{(i)}(\bm{x}_e) - \delta^*}) -1 \right) \right) \leq 0
        \label{eq:msd_cons}
    \end{equation}

\textbf{Latent Space Constraint :} Observe that while the optimizer continuously varies the latent coordinates during optimization, the shapes encoded from the library occupy discrete points in the latent space. In other words, the optimizer is allowed to explore shapes outside the library during optimization. However, at convergence, we expect to select shapes exclusively from the library. To enforce this, we impose a constraint on the distances between the latent coordinates of selected shape instances and those of the library shapes.
With $\mathcal{D}$ being the pairwise distances between the latent coordinates of shape instances and library shapes:

\begin{equation}
    \mathcal{D}_{ji} = ||\bm{z}_*^{(j)} - \bm{z}^{(i)}|| \quad , \; j=1,\ldots,n_L  \; , \; i=1,\ldots,n_F
    \label{eq:dist_latent_part_design_lib}
\end{equation}
We can express the constraint as:

\begin{equation}
    g_l \equiv 
    \underset{i}{\text{max}}( \underset{j}{\text{min}} \; \mathcal{D}_{ji}) \leq 0
    \label{eq:latent_dist_cons}
\end{equation}

To facilitate gradient-based optimization, we use the LogSumExp approximations of the max and min functions.

\textbf{Bound Constraints :} Given that the latent space coordinates $\bm{z} \sim \mathcal{N}(0, 1)$, we constrain $z_{i} \in [-3, 3] \; , \; \forall i$. Further, the translation coordinates $\left(\bar{x}^{(i)}, \bar{y}^{(i)}\right)$ are constrained to lie within the bounding box of the design domain (\cref{fig:variation_transform_parameters}). The orientation parameter is constrained as $0 \leq \theta^{(i)} \leq 2\pi$. Furthermore, the scaling parameter is constrained to $s_{\text{min}} \leq s^{(i)} \leq s_{\text{max}}$. The lower bound $s_{\text{min}}$ ensures the MFS constraint is met. For example, with a library MFS of $40$ nm, imposing an MFS of $60$ nm requires $s_{\text{min}} = 1.5$ (see \cref{sec:results_mfsMsd}). The upper bound $s_{\text{max}}$, the ratio of the diagonal lengths of the design and shape bounding boxes, ensures the shape instance does not exceed the design domain. 

\textbf{Optimization :} Collecting the objective (\cref{eq:objective}), solver (\cref{eq:FDFD_simulation}), and constraint (\cref{eq:msd_cons,eq:latent_dist_cons}) the optimization problem can be expressed as:

\begin{subequations}
	\label{eq:optimization_base_Eqn}
	\begin{align}
		& \underset{\overline {\bm{\Lambda}} = \{ \bm{\Lambda}^{(1)}, \bm{\Lambda}^{(2)}, \ldots \bm{\Lambda}^{(n)} \}} {\text{minimize}} 
           & L(\bar{\bm{\Lambda}}) \label{eq:optimization_base_objective} \\
		& \text{subject to}
		& \bm{K}_{\omega}(\overline {\bm{\Lambda}}) \bm{E} = \bm{J}
            \label{eq:optimization_base_govnEq}\\
            & &  g_s(\overline {\bm{\Lambda}}) \leq 0  \label{eq:optimization_base_separation_Cons} \\
		& &  g_l(\overline {\bm{\Lambda}}) \leq 0
            \label{eq:optimization_base_latentCons} \\
            && \bar{\bm{\Lambda}}_{\text{min}} \leq \bar{\bm{\Lambda}}\leq \bar{\bm{\Lambda}}_{\text{max}} \label{eq:optimization_base_boxCons}
	\end{align}
\end{subequations}

Where $[\bar{\bm{\Lambda}}_{\text{min}}, \bar{\bm{\Lambda}}_{\text{max}}]$ collectively represent the bound constraints on the design variables. The method of moving asymptotes (MMA) \cite{svanberg2007mma}; a gradient-based constrained optimizer is employed to perform the design updates. Finally, to update the design variables $\bar{\bm{\Lambda}}$, we require the gradients of the objectives and constraints with respect to these variables. Leveraging the automatic differentiation (AD) capabilities of JAX \cite{jax2018github}, we automatically derive the sensitivities \cite{chandrasekhar2021auto, ian2020AD}. In practice, this means we only define the forward expressions, and JAX's $\textit{autograd}$ library computes all necessary derivatives with machine precision.

\subsection{Algorithm}
\label{sec:method_algorithm}

Having defined all components of the framework, we now summarize and present the complete algorithm.

As mentioned earlier, we assume a library of shapes has been provided. The procedure begins with training the VAE as detailed in \cref{sec:method_libraryRep_convoImpNN}. The network is trained till sufficiently high representational accuracy is attained. The convolutional encoder is then discarded and the implicit encoder and decoder are retained ($VAE^*$). The network now takes the latent space coordinates and spatial coordinates as input and returns the SDF of the shape at the pertinent latent and spatial coordinates.

The main optimization algorithm is summarized in \cref{alg:trans_opt}. Initially, the design domain is discretized, and the coordinates of the element centers are computed (\ARef**{alg:photos_discretizeDomain}). The design variables are then randomly initialized (\ARef**{alg:photos_init}). We iterate the optimization process until convergence. During each iteration, the latent space coordinates and transformation parameters of the shape instances are obtained from the design variables (\ARef**{alg:photos_extractDesVars}). The transformed coordinates of the design space are computed using the transformation parameters (\ARef**{alg:photos_transformCoordns}). The SDFs of the shape instances are then computed using the latent space coordinates and the transformed coordinates (\ARef**{alg:photos_fwdProp}). The SDFs are then projected to obtain the density fields (\ARef**{alg:photos_sdfToDensity}). The design density field is then computed as a union of the density fields of the shape instances (\ARef**{alg:photos_unionDensity}). We then determine the permittivity field from the density field (\ARef**{alg:photos_densityToPermittivity}). We utilize the permittivities to simulate the optical response of the design (\ARef**{alg:photos_permittivityToElecFeild}). The scattering spectra are computed using the optical response (\ARef**{alg:photos_scatteringSpectra}). Then the objective, minimum separation distance and latent space constraints are computed (\ARef**{alg:photos_objective}, \ARef**{alg:photos_separationConstraint}, \ARef**{alg:photos_latentCons}). The sensitivities are computed in an automated fashion (\ARef**{alg:photos_Sensitivity}) and the design variable is updated using MMA (\ARef**{alg:photos_MMA}). The process continues until the MMA tolerances are met or the iterations exceed a maximum value (\ARef**{alg:photos_iterTillConvg}).

\begin{algorithm*}[]
	\caption{PhoTOS}
	\label{alg:trans_opt}
	\begin{algorithmic}[1]

		\State $\Omega^0_h \rightarrow \bm{x}_e$
            \Comment{elem center coordinates $\bm{x}_e$ of  size $\{n_e \times 2\}$}
            \label{alg:photos_discretizeDomain}

            \State Initialize $\bm{\overline{\Lambda}}^0$
            \Comment{random initialization}
            \label{alg:photos_init}

            \State k = 0

		\Repeat \Comment{optimization loop}

             \State $\overline{\bm{\Lambda}} = \{\bm{z_{1}}, \bm{z_{2}}, \bm{\bar{x}}, \bm{\bar{y}}, \bm{\theta}, \bm{s}\}$
             \label{alg:photos_extractDesVars}

            \State $(\bm{\bar{x}}, \bm{\bar{y}}, \bm{\theta}, \bm{s},  \bm{x}_e ) \rightarrow \tilde{\bm{x}} $
            \Comment{Affine transformations, \cref{sec:method_Transformation}}
            \label{alg:photos_transformCoordns}
            
            \State $(\bm{z}, \tilde{\bm{x}}, VAE^*) \rightarrow \hat{\bm{\phi}}$
            \Comment{shape instance SDFs, \cref{sec:method_libraryRep_convoImpNN} }
            \label{alg:photos_fwdProp}

            \State $\hat{\bm{\phi}} \rightarrow {\bm{\rho}}$
            \Comment{projection of shape SDF to density, \cref{eq:density_projection}}
            \label{alg:photos_sdfToDensity}

            \State $\bm{\rho} \rightarrow \hat{\bm{\rho}}$
            \Comment{design density, \cref{eq:density_union}}
            \label{alg:photos_unionDensity}

            \State $\bm{\rho} \rightarrow \bm{\varepsilon}$
            \Comment{design permittivity, \cref{eq:linear_material} }
            \label{alg:photos_densityToPermittivity}

            \State $ \bm{K}_{\omega}(\bm{\varepsilon}), \bm{J} \rightarrow \bm{E}$
            \Comment{ EM simulation, \cref{eq:FDFD_simulation}}
            \label{alg:photos_permittivityToElecFeild}

            \State $\bm{E} \rightarrow \bm{S}$
            \Comment{ scattering spectra, \cref{sec:method_optimization}}
            \label{alg:photos_scatteringSpectra}

            \State  $(\bm{S}, \bm{S}^{*}) \rightarrow L$
            \Comment{objective, \cref{eq:objective}}
            \label{alg:photos_objective}

            \State $\hat{\bm{\phi}}, MSD \rightarrow g_s$
            \Comment{MSD constraint, \cref{eq:msd_cons}}
            \label{alg:photos_separationConstraint}

            \State $(\bm{z}, \bm{z}_*) \rightarrow g_l$
            \Comment{latent space constraint, \cref{eq:latent_dist_cons}}
            \label{alg:photos_latentCons}

            \State Compute $\nabla_{\bm{\Lambda}}J, \nabla_{\bm{\Lambda}}g_s, \nabla_{\bm{\Lambda}}g_l$
            \Comment{automatic differentiation, \cref{sec:method_optimization}}
            \label{alg:photos_Sensitivity}

            \State MMA$(\bm{\Lambda}^k, L, g_s, g_l, \nabla_{\bm{\Lambda}}L,  \nabla_{\bm{\Lambda}}g_s, \nabla_{\bm{\Lambda}}g_l) \rightarrow \bm{\Lambda}^{(k+1)}$
            \Comment{MMA update step, \cref{sec:method_optimization}}
            \label{alg:photos_MMA}

		\State $\text{k}++$
		
		\Until{ Convergence and k < max\_epoch}
            \label{alg:photos_iterTillConvg}
        \end{algorithmic}
\end{algorithm*}

\section{Numerical Experiments}
\label{sec:results}

In this section, we present several experiments to demonstrate the proposed framework. Without loss of generality, the default parameters for the experiments are set as follows:

\begin{enumerate}
    \item The computational domain consists of a grid of $80 \times 80$ elements, representing a design domain of $1600$ $\times$ $1600 \; \text{nm}^2 $.
    \item A MFS of 40 nm and a MSD of 60 nm is imposed on the design.
    \item The design is populated with 36 shape instances initialized on a $6 \times 6$ grid. The latent coordinates are initialized randomly with a seed of 27.
    \item All components are optimized for their performance within two $10$ nm wavelength bands centered at $\omega = 1270$ nm and $\omega = 1290$ nm (O-band). Excitation source in their fundamental mode is considered.
    \item Designs with an insertion loss above -0.5 dB and a back reflection loss below -20 dB are desired.
    \item The input and output port waveguides have a width of $400$ nm.
    \item Optimization is performed using MMA with a move limit set to $10^{-2}$ for a maximum of 150 iterations. All other default parameters correspond to the version of MMA presented in \cite{svanberg2007mma}.
    \item All experiments are conducted on a MacBook M3 Pro, using the JAX library \cite{jax2018github} in Python. 
\end{enumerate}

\subsection{Convergence}
\label{sec:results_convergence}

\begin{figure}

        \centering
 	\begin{center}
		\includegraphics[scale=0.65,trim={30 10 30 20}]{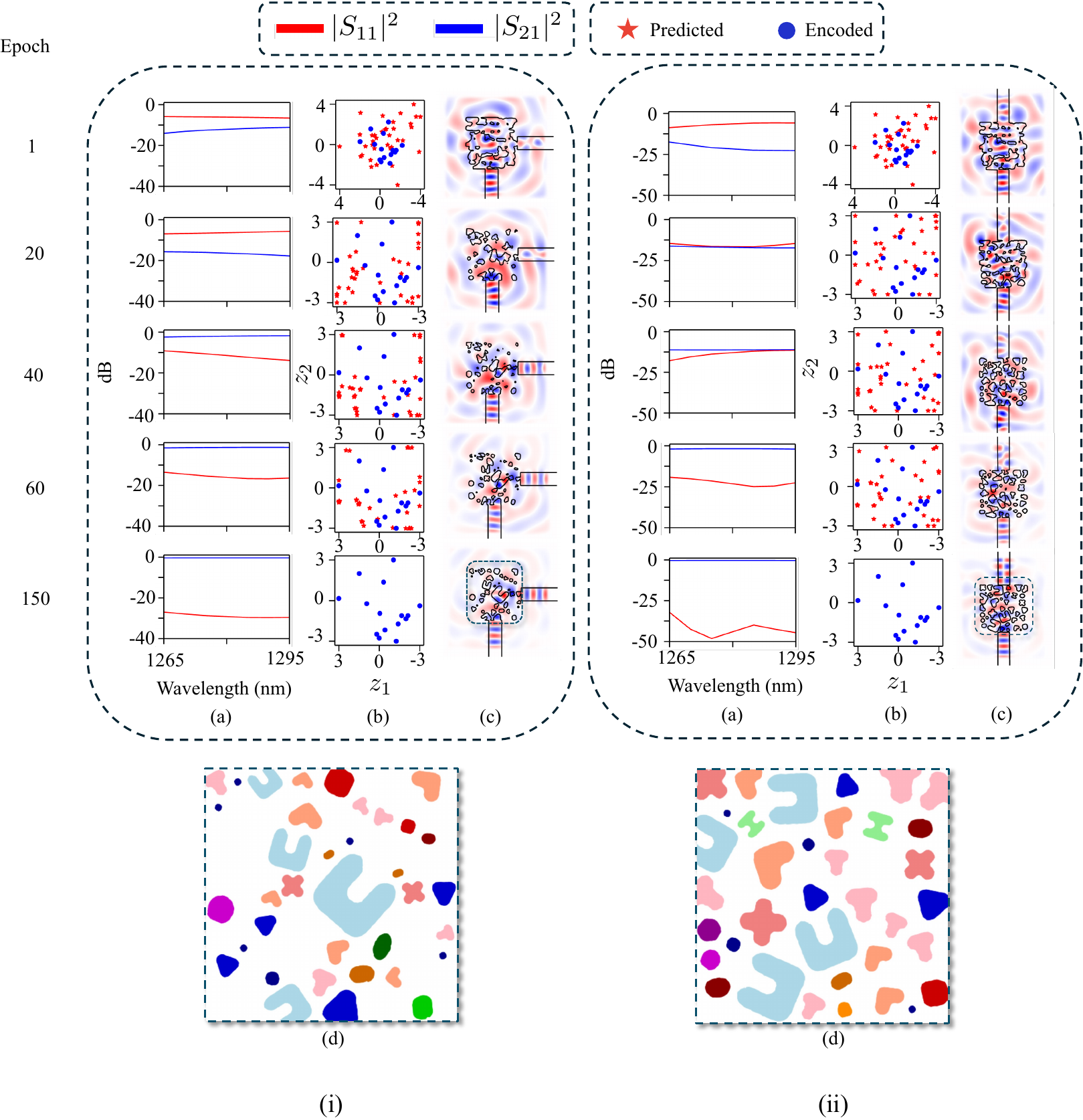}
 		\caption{Evolution of designs:
                    (i) waveguide bend and (ii) mode converter.
                    (a) Scattering spectra.
                    (b) Latent space coordinates of selected and encoded library shapes.
                    (c) Electric field magnitude at wavelength of 1280 nm.
                    (d) Final design.}
        \label{fig:convergence}
	\end{center}
 \end{figure}

In this section, we apply our TO framework to the design of waveguide bends and mode converters. \cref{fig:convergence}(i) depicts the convergence of a waveguide bend. Given an excitation from port 1 (bottom), our objective is to maximize transmission to port 2 (right) while minimizing back reflection. \cref{fig:convergence}(ii) illustrates the convergence of a mode converter design, where the goal is to achieve maximum conversion of the fundamental waveguide mode at port 1 (bottom) to the second-order mode at port 2 (top), with minimal back reflection.

Initially, both designs exhibit poor performance, characterized by low transmission to the output port and significant reflection at the input port. Further, the latent coordinates of the shapes do not align with those of our library shapes. However, subsequent updates substantially modify the design topology by optimizing the transformation parameters and latent coordinates.

By the $150^{th}$ iteration, the target performance specifications are met. Additionally, the separation and latent coordinate constraints are satisfied, with the latent coordinates of the shapes in the design completely overlapping those of our library shapes. In essence, we obtain designs that fulfill the performance requirements while adhering to the MFS and MSD constraints. \cref{fig:convergence}(d) confirms that the final designs consist of well-separated shapes from our library, satisfying the MSD constraint.

\subsection{Feature Size Constraints}
\label{sec:results_mfsMsd}

 \begin{figure}
        \centering
 	\begin{center}
		\includegraphics[scale=0.55]{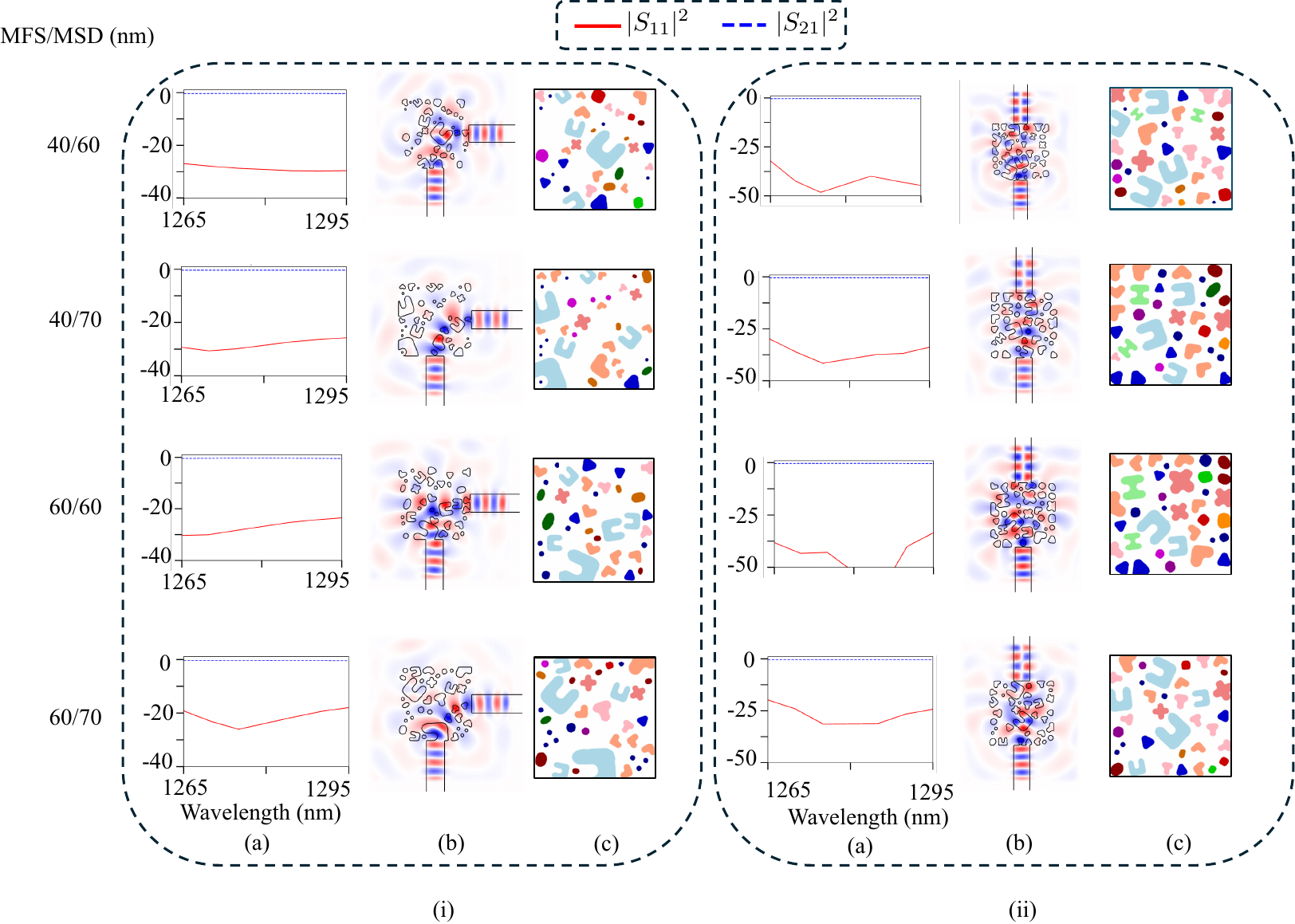}
 		\caption{Performance and design for various minimum feature size (MFS) and minimum separation distance (MSD). (i) Waveguide bend. (ii) Mode converter. (a): Scattering spectra at the input and output port. (b): Electric field magnitude at a wavelength of 1280 nm. (c): Optimized design.}
        \label{fig:constraint_var}
	\end{center}
 \end{figure}
 
A key consideration during the design stage is investigating the impact of the feature sizes; MFS and MSD. While larger feature sizes may enhance fabricability and be necessary for certain fabrication processes, they also restrict the solution space. Recall that the MFS is enforced through a lower bound on the scaling range, and the MSD is imposed via \cref{eq:msd_cons}.

We explore various combinations of MFS/MSD in \cref{fig:constraint_var} for the waveguide bend and mode converter. In all cases, our framework successfully satisfied the imposed size constraints and met performance requirements. However, it is important to note that if the required feature sizes are large, the optimizer may fail to achieve the desired performance due to a lack of feasible points in the solution space. In such scenarios, exploring optimization with larger component dimensions may be necessary.

\subsection{Effect of Initialization}
\label{sec:results_init}

 \begin{figure}[H]
        \centering
 	\begin{center}
		\includegraphics[scale=0.6]{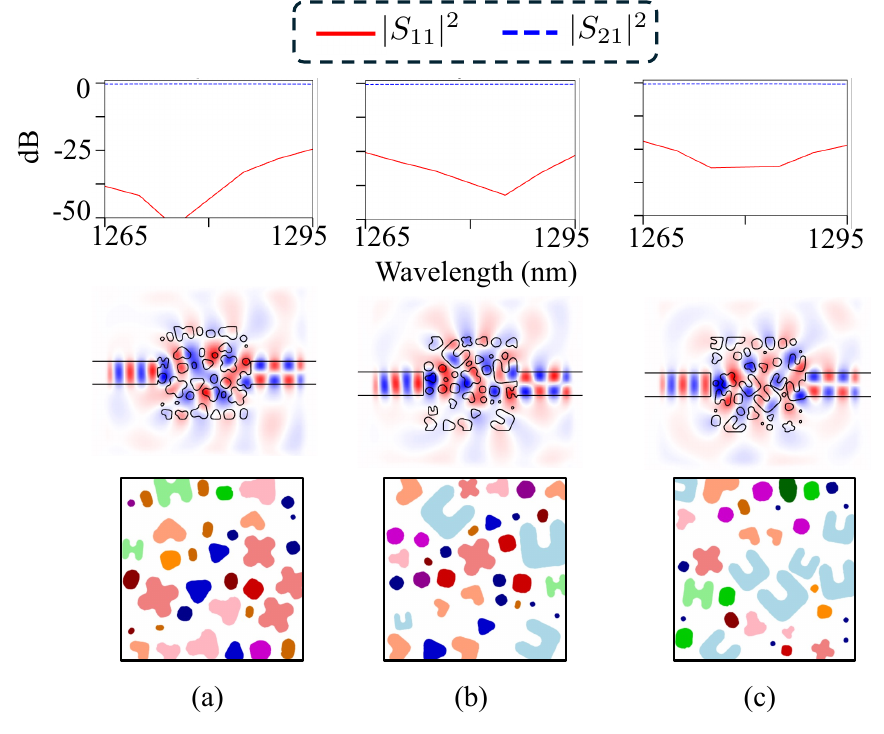}
 		\caption{Optimized performance and design obtained from different random initialization seed of (a) 10 (b) 100 (c) 1000. Top: Scattering parameters at the input and output port. Middle: Electric field magnitude at wavelength of 1280 nm. Bottom: Optimized design.}
        \label{fig:initialization}
	\end{center}
 \end{figure}

In this experiment, we investigate the influence of initial design on the optimization result. Given the non-convex, non-linear nature of our optimization problem, we anticipate that different initializations will lead to distinct local optima. We initialize the optimization with a grid of equi-spaced shapes, setting the scaling to the mean of the MFS and domain diagonal length. The latent coordinates are initialized randomly as $\mathcal{N}(\mu = 0,\sigma = 2)$.

For example, consider \cref{fig:initialization}. Here we optimize mode convertors with varying seeds for the latent coordinates. Observe that while we obtain diverse topologies, the performances are similar; as indicated by their scattering parameters. This suggests that, as expected, the loss landscape is highly non-convex with numerous local solutions. Additionally, it indicates the robustness of our TO framework, capable of discovering designs with the desired performance regardless of initialization.

\section{Conclusion}
\label{sec:conclusion}

In this work, we presented a topology optimization framework for the design of fabricable photonic components. Building upon feature mapping methods, our approach extends their capabilities to encompass multiple generic shapes. Further, we incorporate fabrication constraints, namely minimum feature size and separation in our design. We leveraged a convo-implicit variational autoencoder (VAE) to transform the discrete shape library into a continuous and differentiable latent space, facilitating gradient-based optimization. The effectiveness of our framework was demonstrated through the successful design of waveguide bends and mode converters, highlighting its ability to generate high-performance photonic components that adhere to strict fabrication constraints.

We also identify several avenues for improvement. While the current convo-implicit VAE is trained solely on SDF images, incorporating the Eikonal equation into the training process could enhance accuracy \cite{sitzmann2020SIREN}. Further, the number of shape instances was predetermined in this study. Future work will investigate the inclusion of this number as a design variable, offering greater flexibility. In addition, our library was limited to approximately 15 shapes for illustrative purposes. Expanding the library to include a wider variety of shapes would increase the design space and potentially yield even better results. Furthermore, our current implementation is focused on 2D designs. Extending the framework to 3D designs and validating the fabricated compoents' performance is of significant interest. Finally, we aim to leverage the expanded design space offered by our approach to design more complex photonic components \cite{piggott2015inverse, chang2018inverse, yesilyurt2021efficient}.

\section*{Acknowledgments}
No external funding was used in supporting this work.

\section*{Compliance with ethical standards}
The authors declare that they have no conflict of interest.

\section*{Replication of Results}
The Python code is available at \href{https://github.com/aadityacs/PhoTOS}{github.com/aadityacs/PhoTOS}

\bibliographystyle{unsrt}  
\bibliography{6_references}

\end{document}